\documentclass{article}
\usepackage{graphicx}
\usepackage{amsmath}
\usepackage{amssymb}
\usepackage{mathrsfs}
\usepackage{graphics}
\usepackage{caption}
\usepackage{subcaption}
\usepackage{mathtools}
\usepackage{floatrow}
\usepackage{array}
\usepackage{bbm}
 \usepackage{booktabs}
 \usepackage{float}
 \usepackage[labelfont=bf]{caption}
 \usepackage{varioref}
 \usepackage[T1]{fontenc}
 \usepackage[margin= 1 in]{geometry}
 \usepackage{ucs}
\usepackage[utf8x]{inputenc}

\newcommand{\NN}{\mathcal N}

\newcommand{\XX}{\mathcal X}

\newcommand{\R}{\ensuremath{\mathbb{R}}}

\renewcommand{\P}{\mathbb P}

\newcommand{\vero}[1]{{\color{black}#1}}

\numberwithin{equation}{section}

\begin{document}

\title{Censored pairwise likelihood-based tests for mixing coefficient of spatial max-mixture models{}}

\author{\vero{Abu-Awwad} Abdul-Fattah \and  Maume-Deschamps V\'{e}ronique \and  \vero{Ribereau} Pierre }

\maketitle

\begin{abstract}
	\vero{Max-mixture processes are defined as $Z=\max(aX\/,(1-a)Y)$ with $X$ an asymptotic dependent (AD) process, $Y$ an asymptotic independent (AI) process and $a\in[0\/,1]$. So that, the 
mixing coefficient $a$ may reveal the strength of the AD part present in the max-mixture process.  In this paper we focus on two} tests based on censored pairwise likelihood estimates.  We compare 
their performance through an extensive simulation study. Monte Carlo simulation plays a fundamental tool for asymptotic variance calculations. We apply our tests to daily precipitations from the East 
of Australia. Drawbacks and possible developments are discussed.\\

Keywords: Composite likelihood, Max-stable process; Max-mixture models; Pairwise likelihood; Monte Carlo simulation.

\end{abstract}


\section{Introduction}

The rise of risky environmental events leads to \vero{renewed interest in the} statistical modelling of extremes, for example modelling extreme precipitation is pivotal in flood protection. In the 
last decade, max-stable (MS) models have arised as a common tool for modeling spatial extremes, since they extend the gereralized extreme value (GEV) distribution to the spatial setting, providing a 
consistent multivariate distributions for maxima in arbitrary dimensions.\\

\vero{Max-stable (MS) processes} for spatial data were first constructed using the spectral representation of \cite{de1984spectral}. Subsequent developments have been done on the 
construction of MS process models, \cite{smith1990max}, \cite{schlather2002models}, \cite{kabluchko2009stationary}, and \cite{davison2012geostatistics}.
\\

 Despite many attractive properties of MS models, these processes are restricted \vero{since} only AD or exact indpendence can be modeled. This drawback is constraining when modeling the tail 
behavior of the multivariate distribution of the data, since it is difficult to asses in practice whether a data set should be modeling using \vero{asymptotic dependent (AD) or asymptotic independent 
(AI) models} \cite{thibaud2013threshold} and \cite{davison2013geostatistics}. \vero{In \cite{doi:10.1093/biomet/asr080}, a flexible class of models which may take into account AD and AI is 
proposed. Theses models, called max-mixture (MM), are a mixture of a MS process and an AI process: $Z=\max(aX\/,(1-a)Y)$, with $X$ an AD process, $Y$ an (AI) process and $a\in[0\/,1]$. So that, $a$ 
represents the proportion of AD in the process $Z$ and our purpose is to propose statistical tests on the value of $a$.}\\
\ \\
\textcolor{black}{To the best of our knowledge, only a few papers address testing AI issue for spatial extreme fields.}  \vero{E.g., in \cite{bacro2010testing} a statistical test for \textcolor{black}{AI} of a bivariate maxima vector is proposed and a generalization to spatial context is done. It is based on the  $F$- madogram \cite{cooley2006variograms}.}
Standard likelihood inference \vero{on the parameters of} MS models is not possible in general, since \vero{the full likelihood is not easily computable for MS vector in dimension greater than 
$2$}. Composite likelihood are usefull when the fully specified likelihood is computationally cumbersome, or when a fully specified model is out of 
reach \cite{lindsay1988composite} and \cite{varin2008composite}. Maximum pairwise likelihood estimation for MS models weren first suggested by \cite{padoan2010likelihood} \vero{and is now widely 
used. In particular, in \cite{doi:10.1093/biomet/asr080} and \cite{bacro2016flexible} the parameter inference for MM processes have been driven.} 
\\

\vero{In this paper, we propose a test on the mixing value $a$ of a MM process. This is achieved using pairwise likelihood statistics.} The paper organised as follows. Section 
\ref{sec:spatial_extreme}  reviews the theory of spatial extreme processes MS and MM. The censored pairwise likelihood approach is presented for the statistical inference in Section 
\ref{sec:likelihood}, our testing proposal approach and the main properties are detailed in Section \ref{sec:pair}. In Section \ref{sec:simu} we show, by means of a series of simulation 
studies the performance of our proposed tests. In section \ref{sec:rainfall} we illustrate our testing approach by the analysis of daily precipitation from the East of Australia. Concluding remarks 
and some \vero{perspectives} are addressed in Section \ref{sec:conclusion}.


\section{Spatial extremes modeling}\label{sec:spatial_extreme}
\vero{Throughout our work, $X=\{X(s)\/, s\in \XX\}$, $\XX\subset \R^d$ (generally, $d=2$) is a spatial process will be assumed to be stationary and isotropic.}
 \subsection{Max-stable processes}
   Let $\{X_{k}(s): s \in \mathcal{X} \subset \mathbb{R}^{d} \}$, $k=0,1,2,...,$ are independent replicates of a stochastic process $X$. Then $X$ is a MS 
process if and only if there exist a sequence of continuous functions $\{a_{n}(s)>0\}$ and $\{b_{n}(s)\}$ such that the rescaled process of maxima, $a^{-1}_n(s) \{\max(X_1,X_2,...,X_n)-b_{n}(s)\}$, 
converges in distribution to $X(s)$ \vero{(see \cite{de2006spatial} for more details)}. By this definition, MS \vero{processes offer} a natural choice for modeling spatial extremes.
  \\

\vero{The} univariate extreme value theory, \vero{implies that} the marginal distributions of $X(s)$ are Generalized Extreme value (GEV) distributed, and without  loss of generality the margins can 
transformed to a simple MS process called standard Fr\'{e}chet distribution, $\mathbb {P}(X(s)\leq z)= \exp \{-z^{-1}\}$.
   \\

  Following \cite{de1984spectral} and \cite{schlather2002models}, a \vero{simple} MS process $X(s)$ \vero{has} the following representation
  
  \begin{equation}
  \label{representation}
 X(s)= \max_{k\geq 1} Q_{k}(s)/P_{k}, \qquad s \in \mathcal{X}.
  \end{equation}
  
where $Q_{k}(s)$  are independent replicates of a non-negative stochastic process $Q(s)$ with unit mean at each $s$, and $P_{k}$ are the points of a unit rate Poisson process $(0,\infty)$.
 \\ 

\vero{The joint distribution function of the process $X(s)$ at $K$ locations $s_1\/,\ldots\/, s_K \in\XX$ is given by
\begin{equation} \label{eq:distrib}
\mathbb {P}(X(s_1)\leq z_1\/, \ldots \/, X(s_K)\leq z_K) =  G_{s_1\/,\ldots\/,s_K}(z_{1},...,z_{K})= \exp\{-V_{s_1\/,\ldots\/,s_K}(z_{1},...,z_{K})\}
\end{equation}
where $V_{s_1\/,\ldots\/,s_K}(.)$ is called} the exponent measure. \vero{It} summarises the structure of extremal dependence and satisfies the property of homogenity of order $-1$ and 
$V_{s_1\/,\ldots\/,s_K}(\infty,...,z,...,\infty)=z^{-1}$. 
\vero{It has to be noted that
$$\mathbb {P}\{X(s_1)\leq z,..,X(s_K)\leq z\}=\exp\{-V_{s_1\/,\ldots\/,s_K}(1,...,1)/z\}=\exp\{-1/z\}^{\theta_{s_1\/,\ldots\/,s_K}}\/,$$ 
with $\theta_{s_1\/,\ldots\/,s_K}=V_{s_1\/,\ldots\/,s_K}(1,...,1)$. The coefficient $\theta_{s_1\/,\ldots\/,s_K}$ }  is known as the extremal coefficient. \vero{It} can be seen as a summary of 
extremal dependence with two boundary values, complete dependence $\theta_{s_1\/,\ldots\/,s_K}=1$, and complete independence, $\theta_{s_1\/,\ldots\/,s_K}=K$. In the bivariate case, the AI 
and AD between a pair of random variables $Z_{1}$ and $Z_{2}$, with marginal distributions $F_{1}$ and $F_{2}$, \vero{may be identified by}
\begin{equation}\label{eq:chi}
\chi = \lim_{u \, \to \, 1^-} P(F_{1}(Z_{1} )> u | F_{2}(Z_{2}) > u ))  \/.
\end{equation}
\vero{The cases}  $\chi=0$ and $\chi> 0$ represent AI and AD, respectively, \cite{joe1993parametric}. \vero{This coefficient is related to the pairwise} extremal 
coefficient $\theta$ through the relation $\chi=2-\theta$.
\\

%
Since both dependence functions $\theta$ and $\chi$ are useless \vero{for AI processes}, \cite{coles1999dependence} proposed a new dependence \vero{coefficient which measures the strength of 
dependence for AI processes:}
\begin{equation}
\bar{\chi}=\lim_{u \, \to \, 1^-}\bar{\chi}(u)=\lim_{u \, \to \, 1^-}\frac{2\log P(F(Z(s))>u)}{\log P(F(Z(s))>u, F(Z(s+h))>u)}-1
\end{equation}
AD (respectively AI) is achieved if and if $\bar{\chi}=1$ (resp. $\bar{\chi}<1$).
\\

 Different choices for the process $Q(s)$ in (\ref{representation}) lead to some useful MS models, commonly used choices are the Guassian extreme value process \cite{smith1990max}, the extremal 
Gaussian process \cite{schlather2002models}, the Brown-Resnick process \cite{kabluchko2009stationary}, and the extremal t process \cite{opitz2013extremal}. Below, we list two specific examples of 
MS process models.
 \\
 
 The storm profile model \cite{smith1990max}, \vero{is} defined by taking $Q_{k}(s)= f(s-\vero{P}_{k})$, where $f$ is a Guassian density with covariance matrix $\Sigma \in \mathbb{R}^{2\times 2}$, 
and $\vero{P}_{k}$ is a homogenous Poisson process. The 
 bivariate marginal probability distribution of the Smith model \vero{has} the form

 \begin{eqnarray*}
 \vero{\P(X(s)\leq z_1 \/, X(s+h)\leq z_2)=G_h(z_{1},z_{2})} = \exp \left\{-\left [\frac{1}{z_{1}}\Phi\left( \frac{\gamma(h)}{2}+\frac{1}{\gamma(h)}\log\frac{z_{2}}{z_{1}}\right )+ \right.\right. \\
\left. \left.\frac{1}{z_{2}}\Phi\left( \frac{\gamma(h)}{2}+\frac{1}{\gamma(h)}\log\frac{z_{1}}{z_{2}}\right)\right]\right\}
 \end{eqnarray*}
 
 where $\gamma(h)=\sqrt{h^{T}\Sigma^{-1}h}$, in which $h$ is the seperation vector between the two locations,$\ s_{1}$, and $\ s_{2}$, and $\Phi$ is the standard normal distribution. The pairwise 
extremal coefficient is $\vero{\theta_{s\/,s+h}=}\theta(h)=2\Phi\{\gamma(h)/2\}$. 
 \\

 The Truncated Extremal Gaussian (TEG) model  \vero{is originally due to \cite{schlather2002models}, it is defined using}
  
  \begin{equation}
 Q_{k}(s)= c \max (0,\varepsilon_{k}(s))\mathbbm{1}_{\mathcal{A}_{k}}{(s-\vero{P}_{k})}
  \end{equation} 
  
  where $\varepsilon_{k}(s)$ are independent replicates of a stationary Gaussian process $\varepsilon=\{\varepsilon(s),s \in \mathcal{X}\}$ with zero mean, unit variance and correlation function 
$\rho(.)$. $\mathbbm{1}_{\mathcal{A}}$ is the indicator function of a compact random set $\mathcal{A} \subset \mathcal{X}$, $\mathcal{A}_{k}$ are indepndent replicates of $\mathcal{A}$ and 
$\vero{P}_{k}$ are points of Poisson process with a unit rate on $\mathcal{X}$. The constant $c$ is chosen to satisfy the constraint $\mathbb{E}\{Q_k(s)\}=1$.\\

The bivariate marginal probability distribution of TEG model in the stationary case has the form

\begin{eqnarray}
\label{truncated}
\vero{\P(X(s)\leq z_1 \/, X(s+h)\leq z_2)=G_h(z_{1},z_{2})}=\exp \left\{ -\left( \frac{1}{z_1}+ 
 \frac{1}{z_2}\right) \right. \cdot \\
 \left. \left[ 1- \frac{\alpha(h)}{2}\left(1- \sqrt{1- \frac{2 (\rho(h)+1)z_{1} z_{2}} {(z_{1}+z_{2})^{2}}}\right)\right] \right\}\nonumber
\end{eqnarray}
where $h$ is the spatial lag, and $\alpha(h)\vero{=} (1- h/2r) \mathbbm{1}_{[0,2r]}$ \vero{if}  $\mathcal{A}$ \vero{is} a disk of fixed radius $r$. The pairwise extremal coefficient 
$\theta(h)= 2 - \alpha({h}) \{1-\sqrt{(1-\rho(h))/2}\}$. TEG model were fitted by \cite{davison2012geostatistics} to extreme temperature data.\\

 The stochastic process $X(s)$ defined in equation(\ref{representation}) has the bivariate density function
 \begin{equation*}
 g(z_1,z_2)=\left\{\frac{\partial}{\partial z_1}V_{X}(z_1,z_2)\frac{\partial}{\partial z_2}V_{X}(z_1,z_2)- \frac{\partial^2}{\partial z_1 z_2}V_{X}(z_1,z_2)\right\}\exp\{-V_{X}(z_1,z_2)\}
 \end{equation*}
 where $V_{X} $ is the exponent measure of the MS process $X(s)$.
\subsection{Hybrid models of spatial extremal dependence}

Although MS models seems to be suitable for modeling extremely high threshold exceedances, AI models may show a better fit at finite thresholds. Since \vero{it may be difficult or impossible } 
in practice to decide whether a dataset should be modeled \vero{using} AD or AI, \cite{doi:10.1093/biomet/asr080} have introduced \vero{an} hybrid spatial dependence model, which \vero{is} able to 
 capture both AD and AI.\\
Let $ X(s), s \in \mathcal{X} $, be a stationary simple MS process, and $Y(s), s \in \mathcal{X} $, be a stationary AI process with unit Fr\'{e}chet margins \vero{(see below for the construction of 
such a process)}. Then for a mixture proportion $\ a \in [0,1]$, \vero{a spatial max-mixture (MM) process is constructed:}
\begin{equation}
\label{max-mixture}
Z(s) =\max\{\ a X(s),(1-\ a)Y(s)\}
\end{equation}

The bivariate distribution function for a pair of sites is \vero{straightforwardly} obtained by the independence between \vero{the processes $X$ and $Y$}

 \begin{equation}
 \label{Biv MM}
\vero{G_h}(z_1,z_2)=\vero{G_h^X}\left(\frac{z_1}{\ a}, \frac{z_2}{\ a}\right) \vero{G_h^Y}\left(\frac{z_1}{1-\ a},\frac{z_2}{1-\ a}\right)
\end{equation}
\vero{where $G_h^X$ (resp. $G_h^Y$) is the bivariate distribution function of $X$, with space lag $h$ (resp. of $Y$). \\
AI processes with unit Fréchet marginal distributions can easily be constructed (see \cite{doi:10.1093/biomet/asr080}). Consider} $Y(s) = -1/\log(\varPhi(Y'(s)))$,  where $\{Y'(s), s \in 
\mathcal{X}\}$ is a Gaussian process. \vero{Then, $Y$ is an AI process with unit Fréchet marginal distributions. Another} class of AI processes \vero{called} inverted max-stable (IMS) 
processes are defined \vero{using a simple}  MS process $\ X(s)$, \vero{let}

	\begin{equation}
	\label{IMS}
		Y(s) = -1/\log\{{1-\exp{[-X(s)^{-1}]}}\}
	\end{equation}
	
With this construction, \vero{any} MS process may be transformed to provide \vero{an} AI counterpart. Bivariate distribution function and density of the margins of $Y(s)$ are respectively
\vero{\begin{equation}
G_h^Y(z_1,z_2)= -1 +\exp(-z^{-1}_{1})+\exp(-z^{-1}_{2})+\exp\{-V_h^Y{[\omega(z_{1}),\omega(z_{2})]}\}
\end{equation}

\begin{align*}
g_h^Y(z_1,z_2)&=\left\{\frac{\partial}{\partial z_1}V_h^Y[\omega(z_1),\omega(z_2)]\frac{\partial}{\partial z_2}V_h^Y[\omega(z_1),\omega(z_2)]- \frac{\partial^2}{\partial z_1 
z_2}V_h^Y[\omega(z_1),\omega(z_2)]\right\}\\
& \times \exp\{-V_h^Y[\omega(z_1),\omega(z_2)]\} \frac{\omega^{2}(z_1)\omega^{2}(z_2)\exp\{-(z_{1}^{-1}+z_{2}^{-1})\} }{z_{1}^{2}z_{1}^{2}\{1-\exp(-z_{1}^{-1}\}\{1-\exp(-z_{1}^{-1}\}}
\end{align*}
where $V_h^Y$ is the exponent measure of the bivariate distribution $(Y(s),Y(s+h))$, and $\omega(z)=-1/\log\{{1-\exp{[-z^{-1}]}}\}$.\\
}
Thus, in the case where $X(s)$ is a MS process and $Y(s)$ is a IMS process, the distribution function in (\ref{Biv MM}) has the form
 \begin{equation}
 \exp\{- a V_{X}(z_1,z_2)\}\{-1+ \exp[-(1-a)/z_1]+ \exp[-(1-a)/z_2]+\exp[-V_{Y}[\omega({1-a}/{z_1}),\omega(1-a/z_2)]\}
\end{equation}

\cite{bacro2016flexible} analyzed daily rainfall data in the east of Australia with a class of different models (MS, AI, and MM), and showed that MM models has the merit to overcome the limits of MS 
models in which only AD or exact \vero{independence} can be modeled.

\section{Inference for MM processes: censored pairwise likelihood approach}\label{sec:likelihood}
\vero{In order to propose a testing procedure on} on the mixing coefficient \vero{a} of MM processes \vero{defined by equation} (\ref{max-mixture}), \vero{we shall} use the composite 
likelihood.\\

\vero{The} composite likelihood technique \cite{lindsay1988composite} \vero{is} a general method of inference for dealing with large datasets \vero{and/or miss-specified models}. A composite 
likelihood consists of a combination of valid likelihood objects usually related to small subsets of data and defined as
\begin{equation}
\label{composite}
\mathcal{L}_{C} = \prod_{k=1}^{K} g(z\in \mathcal{B}_{k};\varphi)^{\omega_{k}}
\end{equation}
where $\{g(z\in \mathcal{B}_{k};\varphi), z \in \mathcal{Z} \subseteq \mathbb{R}^{K}, \varphi \in \varPhi \subseteq \mathbb{R}^{d}\}$ is a paremetric statistical model, $\{\mathcal{B}_{k}; k=1,...,K\}$ is a set of marginal or conditional events, $\{\omega_{k}; k=1,...,K\}$ is a set of suitable weights, if the weights are all equal they may be ignored, non-equal weights may be used to improve the statistical performance in certain cases. The associated composite log-likelihood is $c\ell(\varphi; y)=log\mathcal{L}_{C}(\varphi; y)$.
\\

 In the spatial setting, the definition of a pairwise log-likelihood is derived from (\ref{composite}) by taking $\mathcal{B}_{k}=\{ z_{k}(s_{j}),z_{k}(s_{j'})\}$ as the set of bivariate subvectors of $z$ taken over all $K(K-1)/2$ disitinct sites pairs $j$ and $j'$. Thus the weighted pairwise log-likelihood is given by

\begin{equation}
p\ell(\vartheta;z)= \sum_{k=1}^{M}\sum_{j=1}^{K-1}\sum_{j'=i+1}^{K} \omega_{jj'}\log     \mathcal{L}(z_{k}(s_{j}),z_{k}(s_{j'});\vartheta))
\end{equation}

where $z$ are the data \vero{available} on the whole region, $z_{k}(s_{j})$ is the $k-$th observation of the $j-$th site, and $\mathcal{L}(z_{k}(s_{j}),z_{k}(s_{j'});\vartheta))$ is the likelihood 
function based on observations at locations $j$, and $j'$, $\omega_{jj'}$ are non negative weights specifying the contributions of each pair. A simple weighting choice is to let $\omega_{jj'}= 
\mathbbm{1}_{\{ h \leq \delta\}}$,
where $h$ is the pariwise distance, and $ \delta$ is a specified value.
\\

 Inference using pairwise likelihood methods \vero{is} computationaly expensive, \vero{since} with $K$ sites there are $\binom{K}{2}$ pairs to include. This methodology has been used by 
\cite{padoan2010likelihood}, \cite{davison2012geostatistics} and \cite{davison2012statistical} for infernence on MS processes.
\\
\\
Different inference \vero{approaches} based on a censored threshold-based log-pairwise likelihood have been 
used by several researchers \cite{ledford1996statistics}, \cite{doi:10.1093/biomet/asr080}, \cite{bacro2016flexible} and \cite{huser2014space}. Where the censored pairwise contributions 
$\mathcal{L}_{u}(z_{k}(s_{j}),z_{k}(s_{j'});\vartheta))$ take the forms

\begin{equation}
\mathcal{L}_{u}(z_{k}(s_{j}),z_{j}(s_{j'});\vartheta)=
\begin{cases}
G(u,u;\vartheta), & \text{if max}\ (z_{k}(s_{j}),z_{k}(s_{j'})) \leq u \\
\frac{\partial}{z_{k}(s_{j})} G(z_{k}(s_{j}),u;\vartheta), & \text{if}\ z_{k}(s_{j})>u,z_{k}(s_{j'}) \leq u \\
\frac{\partial}{z_{k}(s_{j'})} G(u,z_{k}(s_{j'});\vartheta), & \text{if}\ z_{k}(s_{j})\leq u,z_{k}(s_{j'}) > u \\
g(z_{k}(s_{j}),z_{k}(s_{j'});\vartheta), & \text{if max}\ (z_{k}(s_{j}),z_{k}(s_{j'})) \geq u \\
\end{cases}
\end{equation}

\begin{equation}
\label{pairwise1}
\mathcal{L}_{u}(z_{k}(s_{j}),z_{j}(s_{j'});\vartheta)=
\begin{cases}
G(u,u;\vartheta), & \text{if max}\ (z_{k}(s_{j}),z_{k}(s_{j'})) \leq u \\
g(z_{k}(s_{j}),z_{k}(s_{j'});\vartheta), & \text{if max}\ (z_{k}(s_{j}),z_{k}(s_{j'})) \geq u \\
\end{cases}
\end{equation}
where $u \in \mathbb{R}$ is a high threshold, $G(.,.)$ is given in equation (\ref{Biv MM}) and $g(.,.)$ \vero{is the bivariate density function, i.e. $g(x\/,y) =\frac{\partial^2 
G(x\/,y)}{\partial x \partial y}$. }

\section{Pairwise likelihood  statistics for testing $H_{0}:a =a_{0}$ versus  $H_{1}:a \neq a_{0} $}\label{sec:pair}

We assumed the parameters of a MM model $  \vartheta \in \mathbb{R}^{d}$ is partitioned as $\vartheta= ( a,\eta)$, where $a$ is a one-dimensional parameter of interest that denotes the mixing coefficient for a MM model and $\eta$ is a $\ q \times 1$ nuisance parameter, $ d = 1 + q$.\\

\vero{Our purpose is} to test the hypothesis $H_{0}:a =a_{0}$ versus  $H_{1}:a \neq a_{0} $, for some specified value $a_{0} \in$[0,1]. Let ${\hat{\vartheta}}_{p\ell}=( \hat{a},\hat{\eta})$ denotes 
the unrestricted maximum pairwise likelihood \vero{estimator} and ${\hat{\vartheta^{*}}}_{p\ell}=(a_{0},\hat{\eta}(a_{0}))$, $\hat{\eta}(a_{0})$ denotes the constrained maximum pairwise likelihood 
\vero{estimator} of ${\eta}$ for a fixed $a=a_{0}$. The pairwise maximum likelihood estimator ${\hat{\vartheta}_{p\ell}}$ is asymptotically normally distributed:

\begin{equation}
\sqrt{M}({\hat{\vartheta}_{p\ell}}-\vartheta)\xrightarrow{\mathscr{D}} \vero{\NN}_{d}\{0,\mathcal{G}^{-1}(\vartheta)\}
\end{equation}

where $\vero{\NN}_{d}\{\mu,\Sigma\}$ is the $d$-dimensional normal distribution with mean $\mu$ and variance $\Sigma$, and $\mathcal{G}(\vartheta)$ denotes the Godambe information matrix:
$$\mathcal{G}^{-1}(\vartheta)=\mathcal{H}^{-1}(\vartheta)\mathcal{J}(\vartheta)\mathcal{H}^{-1}(\vartheta)\/,$$ 
where  $\mathcal{H}(\vartheta)= \mathbb{E}\{- \nabla^{2}{p\ell(\vartheta)\}}$ is called 
the \vero{sensitivity} matrix and $\mathcal{J}(\vartheta)=\mathbb{V}ar\{ \nabla{p\ell(\vartheta)\}}=\mathbb{E}\{\nabla{p\ell(\vartheta)}\nabla^{T}{p\ell(\vartheta)}\}$ is called the variability 
matrix. For more details see \cite{kent1982robust}, \cite{lindsay1988composite}, and \cite{varin2011overview}.
\\

In this paper we propose the following two statistics exploiting the pairwise maximum likelihood as an inferential tool. Our objective to facilitate the modeling of the spatial data by a random 
field with appropriate extremal behaviour. The $Z$-test statistic which is  \vero{straightforwardly} derived \vero{from the} central limit theorem \vero{for maximum composite likelihood estimators}.

\begin{equation}
	Z=\frac{\hat{a} - a} {\sqrt{\mathcal{G}^{aa}(\hat{\vartheta}_{p\ell})}}\xrightarrow{\mathscr{D}} N\{0,1\}
\end{equation}
where $\mathcal{G}^{aa}(\hat{\vartheta}_{p\ell})$ denotes a $1\times1$ submatrix of the inverse of $\mathcal{G}(\hat{\vartheta}_{p\ell})$ pertaining to $a$. While the pairwise likelihood ratio statistic ($LR$) with a nonstandard asymptotic chi-squared distribution is given by

\begin{equation}
LR=2 \{p\ell({\hat{\vartheta}}_{p\ell})-p\ell({\hat{\vartheta^{*}}}_{p\ell}\} \xrightarrow{\mathscr{D}} \lambda \chi^{2}_{1} 
\end{equation}

	where $\lambda = \{\mathcal{H}^{aa}(\hat{\vartheta^{*}}_{p\ell})\}^{-1}\mathcal{G}^{aa}(\hat{\vartheta^{*}}_{p\ell})$, $\mathcal{H}^{aa}(\hat{\vartheta^{*}}_{p\ell})$ and $\mathcal{G}^{aa}(\hat{\vartheta^{*}}_{p\ell})$ are respectively $1\times1$ submatrices of the inverse of $\mathcal{H}{(\hat{\vartheta^{*}}_{p\ell})}$ and  $\mathcal{G}{(\hat{\vartheta^{*}}_{p\ell})}$ pertaining to $a$ and $\chi^{2}_{1}$ is a chi square random variable with one degree of freedom. 
	\\
	
The asymptotic distribution of the $LR$ statistic \vero{ has been studied in \cite{kent1982robust} in a more general context  
(specifically when the dimension of the parameter of interest may be greater than $1$).} Different \vero{kind} of adjustments were proposed to recover an asymptotic chi square distribution in 
\cite{chandler2007inference}, \cite{geys1999pseudolikelihood}, \cite{kuonen1999miscellanea}, \cite{pace2011adjusting} and \cite{rotnitzky1990hypothesis}. \\

Standard errors and critical values for the tests require the estimation of the Godambe matrix and its components, and since analytical expressions for 
$\mathcal{H}{(\vartheta)}$ and $\mathcal{J }{(\vartheta)}$ are difficult to obtain in \vero{mostly} realistic applications. \vero{They are usually estimated} by means of a Monte Carlo 
\vero{simulations:}
$$\hat{\mathcal {H}}^{S}(\vartheta)= - M^{-1}\sum_{k=1}^{M}{ \nabla^{2}p\ell(\hat\vartheta_{p\ell}; z^{k}})$$
$$\hat{\mathcal{J}}^{S}(\vartheta)= M^{-1} \sum_{k=1}^{M}\nabla{p\ell(\hat\vartheta_{p\ell}; z^{k}}) \nabla^{T}{p\ell(\hat\vartheta_{p\ell}; z^{k}})\/,$$
where $z^{k}$, $k=1,...,M$, denote the $k$th datasets simulated from the fitted model. The results obtained by \cite{cattelan2016empirical} for testing the paremeters of equicorrelated multivariate 
normal model showed that the coverage of the statistics based on Monte Carlo simulation are almost indentical to those of statistics based on analytically computed quantities.\\

\vero{Testing} at boundary points \vero{$\ a=1$ (AD)} or $\ a=0$ (AI) is non standard since there are additional nuisance parameters which \vero{are} present only under the alternative 
\cite{davies1977hypothesis} and \cite{davies1987hypothesis}, we apply our $LR$ test at some points \vero{close} to the boundaries, i.e. $a_{0} = 0.01$ or $a_{0}=0.99$.\\


\section{Simulation study}\label{sec:simu}

We \vero{have performed several}  of simulation studies \vero{in order} to investigate the performance of our \vero{testing procedures. We simulated} from a MM model (\ref{max-mixture}) in which 
$X$ is a TEG process (\ref{truncated}) with $\mathcal{A}_{X}$ a disk of fixed radius $r_{X}$. The AI process $Y$ is an inverse TEG process with $\mathcal{A}_{Y}$ a disk of fixed radius $r_{Y}$. 
For simplicity, we assume that the correlation functions \vero{of these} two processes are exponential, with range parameters $\phi_{X} ,\phi_{Y}>0$ respectively. \vero{Our} purpose \vero{is} to 
test $H_{0}:a =a_{0}$ versus  $H_{1}:a \neq a_{0} $, $a_{0}$ varies from 0.01 to 0.99 by steps of 0.01.\\

The censored pairwise likelihood approach (\ref{pairwise1}) is used for estimation, where the threshold $u$ is taken corresponding to the 0.9 empirical quantile at each site, and equal weights are 
\vero{considered}. \textcolor{black}{To reduce computational burden the pairwise likelihood function has been coded in C; the optimization has been parallelized on 20 cores using the R library parallel and carried out using the Nelder–Mead optimization routine in R.}\\

\begin{figure} [H]
	
	\includegraphics[width=0.99\linewidth, height=10cm]{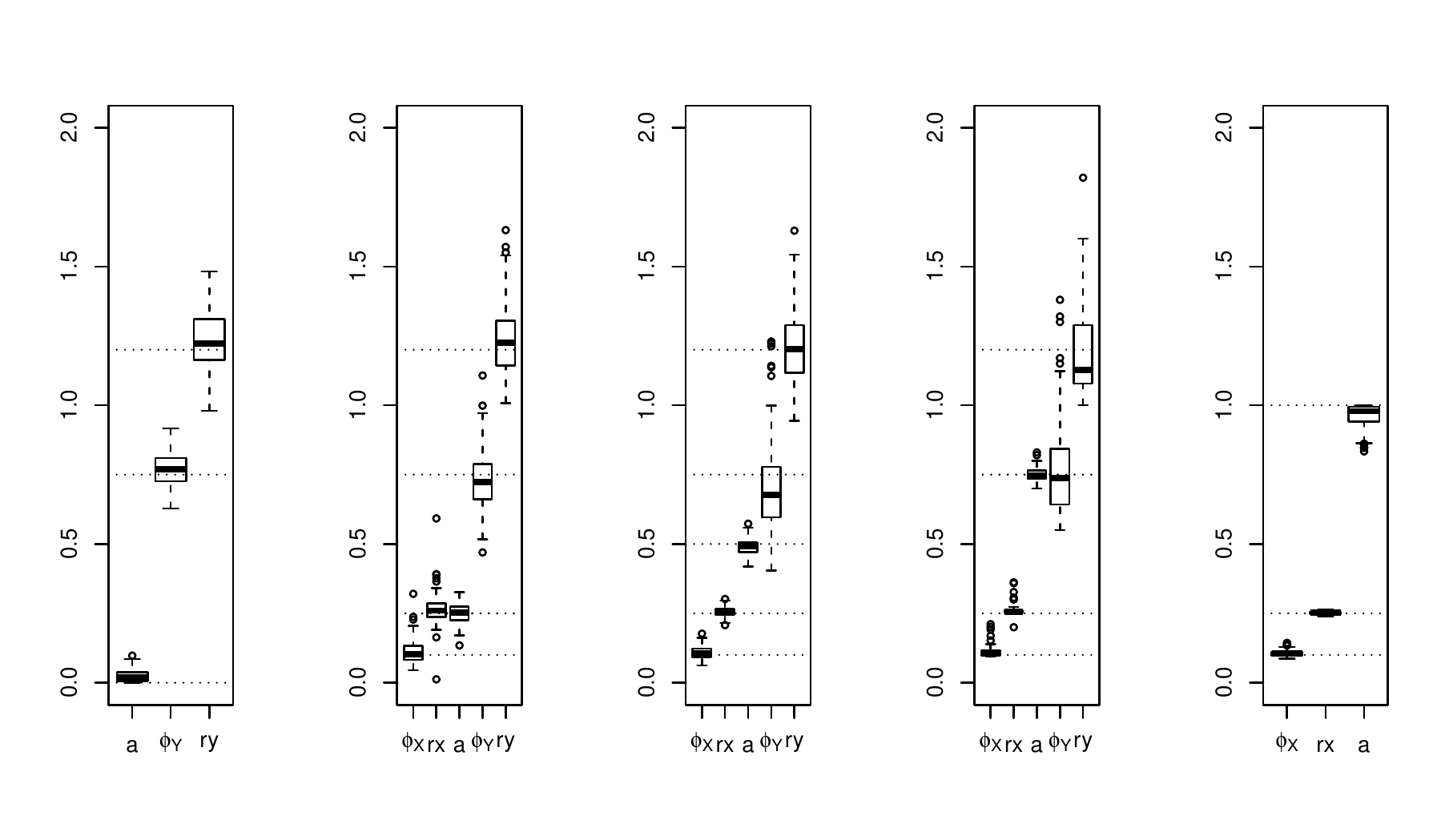}
	
	\caption{Boxplots of censored pairwise likelihood estimates based on 150 simulation replicates from 1000 independent copies of MM model with mixing coefficient (from left to right: $a=0, a=0.25, a=0.50, a=0.75$ and $a=1$) and $\phi_{X}= 0.10,r_{X}=0.25, \phi_{Y}= 0.75, r_{Y}=1.20$. Results for the unidentifiable parameters are not reported. }
	
	\label{boxplots}
	
\end{figure}

\vero{We have done $J=150$ replications of $N=1.000$ independent copies of the considered MM process on $50$ locations uniformaly chosen in the the square $[0, 2]\times [0,2]$. The Boxplots of the 
estimated parameters on the $J$ samples are displayed in Figure \ref{boxplots}. The parameters used are}  $\phi_{X}= 0.10,r_{X}=0.25, \phi_{Y}= 0.75, r_{Y}=1.20$ and different mixing coefficient $ a 
\in \{0,0.25,0,50,0.75$ and $ 1 \}$ \vero{are considered.} \textcolor{black}{Generally, the parameters are well estimated.}


\begin{figure} [H]
	
	\includegraphics[width=0.9\linewidth, height=10cm]{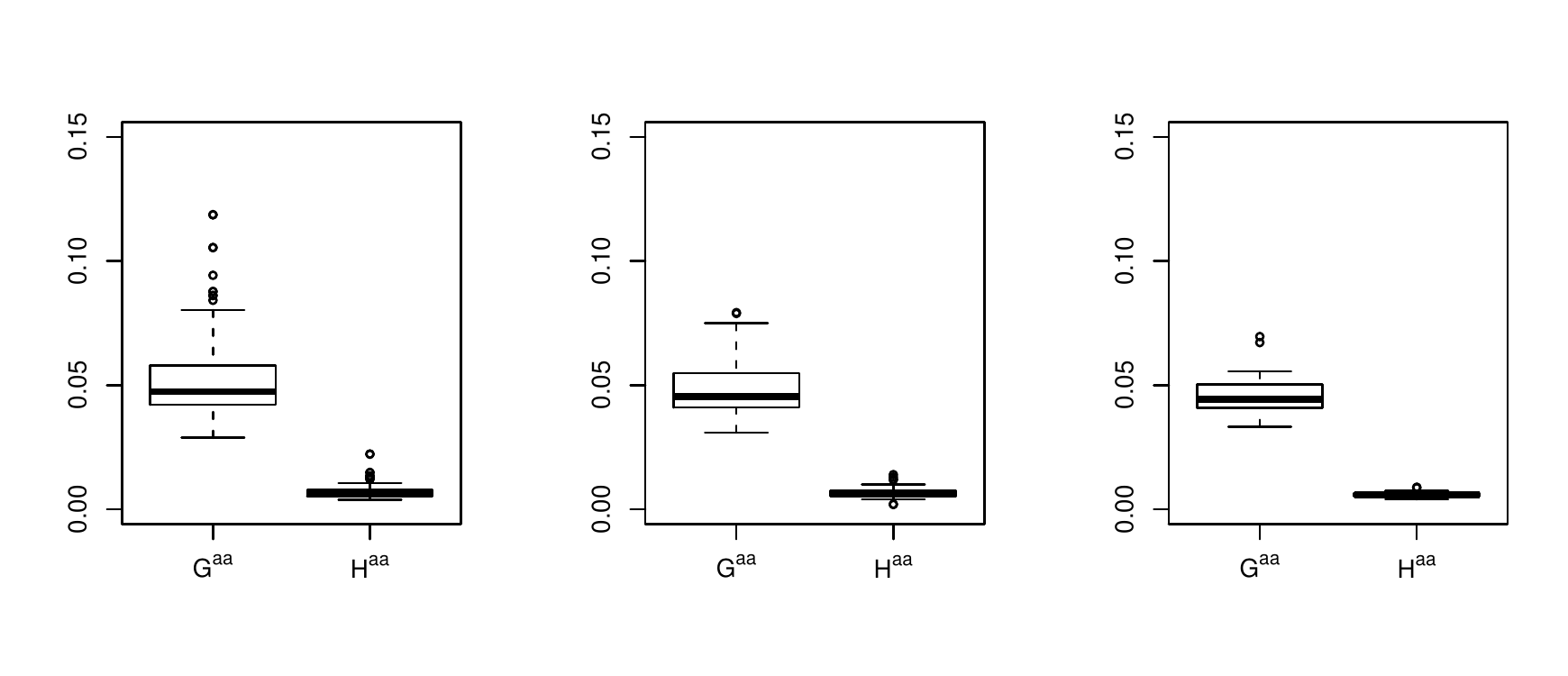}
	
	\caption{Boxplots for calculated matrices elements by Monte Carlo procedure based on 200 simulation replicates from 1000 independent copies of MM model with  parameters $\{\phi_{X}= 0.10, r_{X}=0.25, \phi_{Y}= 0.75, r_{Y}=1.2, a=0.5\}$. (from left to right: $M=1000, M=1250, $ and $ M=1500$ ). Results for required elements  are reported. }
	
	\label{boxplot2}
	
\end{figure}

In order to obtain accurate estimates of $ \mathcal{H}(\vartheta)$ and $\mathcal{J}(\vartheta)$ in the Monte Carlo procedure, we perform an exploratory study with 200 simulation replicates based on 
MM model described above with parameters $\beta_{MM}=\{\phi_{X}= 0.10, r_{X}=0.25, \phi_{Y}= 0.75, r_{Y}=1.2, a=0.5 \}$. For each replication, we randomly generate $K = 50$ locations uniformaly 
in the square [0, 2] $\times $[0, 2] with 1000 independent observations at each sampled location. Then, we simulate data from the fitted model with $M \in \{1000, 1250, 1500\}$ independent 
\vero{simulations}  at the sampled $K$ locations. 
\vero{In Figure \ref{boxplot2}, we present the boxplots for $G^{aa}$ and $H^{aa}$. The results give a justification to use $M=1500$ as a compromise between accuracy and computation time.}

\begin{figure} [H]
	
	\includegraphics[width=0.99\linewidth, height=12cm]{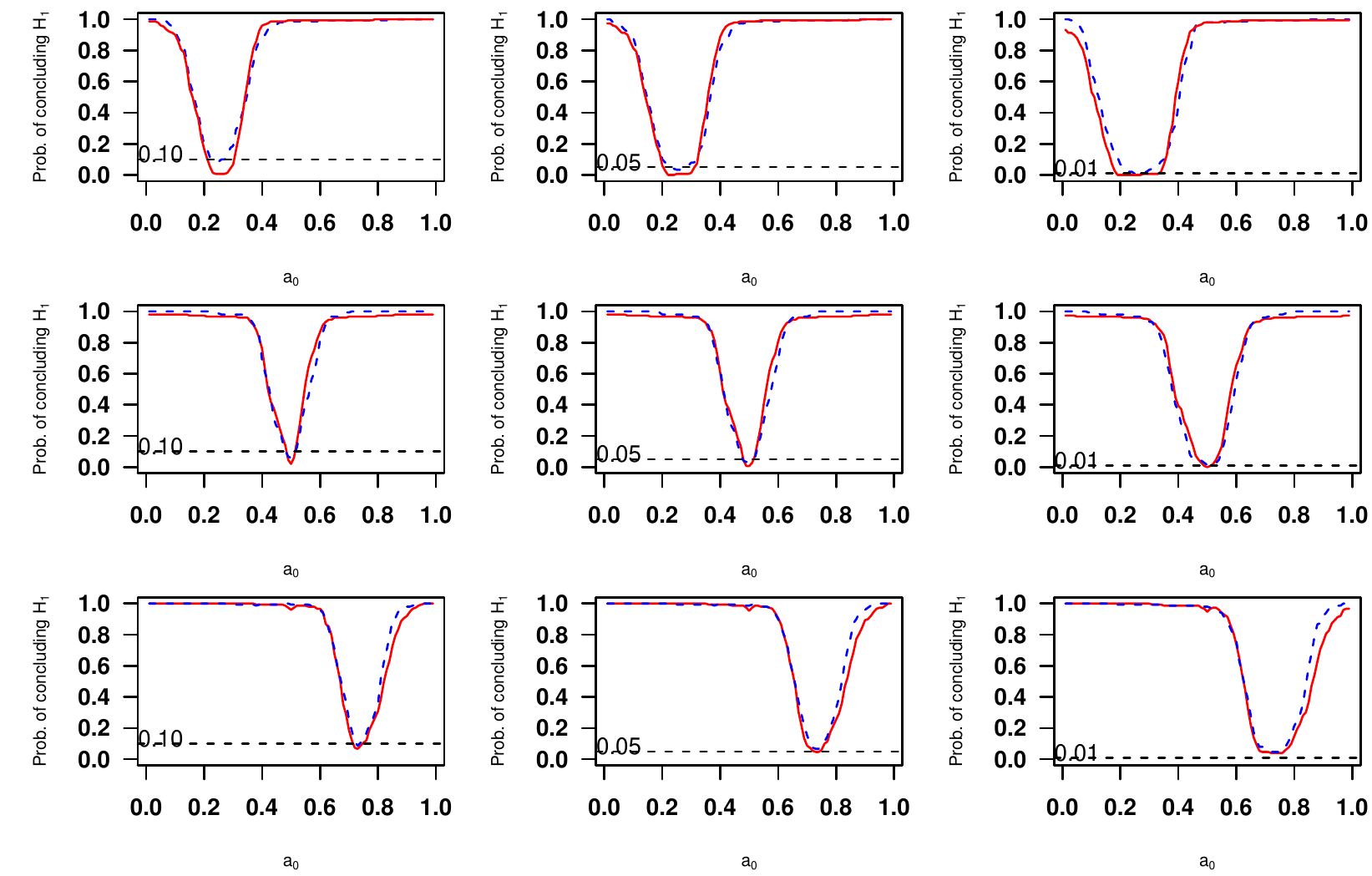}
	
	\caption{ Empirical probabilities of concluding $H_{1}$ based on 150 replicates simulation study of three MM models with $\ X$=TEG, $\beta_{MS}=\{\phi_{X}= 0.10, r_{X}=0.25\}$, and $Y$= 
Inverted TEG processes,  $\beta_{IMS}=\{\phi_{Y}= 0.75, r_{Y}=1.2\}$, the mixing coefficients (top row: $a=0.25$, middle row: $a=0.50$ and bottom row: $a=0.75$). 50 locations were randomly 
and \vero{uniformly} \vero{generated} in the square $[0, 2]^{2}$, and 1000 \vero{realisations} were simulated at the sampled locations. (Red solid line: $Z$ test and blue dashed line: $LR$ test).}
	
	\label{test1}
	
\end{figure}

\vero{Figure \ref{test1} and Tables \ref{table1}, \ref{table2} and \ref{table3} report} a summary of comparison results between the  two statistics $Z$ and $LR$ in terms of empirical 
probabilities for concluding $H_{1}$ in testing hypothesis $H_{0}: a=a_{0}$ against $H_{1}: a \neq a_{0}$ based on 150 simulation replicates from 1000 independent copies simulated at 50 randomly and 
uniformaly sampled locations in the square $[0, 2]^{2}$ from three MM models with parameters $\{\phi_{X}= 0.10, r_{X}=0.25, \phi_{Y}= 0.75, r_{Y}=1.2\}$, according to different values \vero{of} the 
mixing parameter $a \in \{0.25, 0,50, 0.75\}$. Decisions obtained at three significance levels $\alpha \in\{0.01,\vero{0.05},0.10\}$.\\

Despite the poor performance which may be expected at the region around to the true mixing coefficient $a$. The results in terms of empirical probabilities for concluding $H_{1}$ of the two statistic 
\vero{show} a reasonable performance. The probability of making a correct decision becomes higher as we become very close or move far from the model true value.     \vero{We also remark that the 
performances of the two tests are very similar except for $a_0=0.25$, $\alpha=0.1$ where the $Z$ test presents an unexpected over sensitivity.}

\newcommand{\ra}[1]{\renewcommand{\arraystretch}{#1}}


\begin{table}[H]\centering

	\begin{tabular}{@{}cccccccccccc@{}}

		\bottomrule
		$a_{0}$ & 0.05 & 0.10 &0.15 & 0.20&0.25& 0.30&0.35 & 0.40 & 0.50&0.80  \\
		
	\bottomrule

$\alpha$=0.01\\
	\bottomrule
		
		$LR$ & 0.980 & 0.647& 0.293 & 0.067 & 0.007 &0.033  &0.100& 0.440&0.980&0.993\\
		$Z$ & 0.880&0.533 & 0.187& 0.000& 0.000& 0.007& 0.067&0.600 &0.980&0.993\\
		
	\end{tabular}

	\begin{tabular}{@{}cccccrrrcrrr@{}}
		
		\bottomrule
		$\alpha$=0.05\\
		\bottomrule
		$LR$ &0.973  &0.880 &0.480  &0.100  & 0.033 & 0.080 &0.293& 0.813&0.980&0.993\\
		$Z$ & 0.993&0.813 &0.440 &0.060 & 0.007 &0.013 & 0.393& 0.880&0.987&0.993 \\

	\end{tabular}

	\begin{tabular}{@{}cccccrrrcrrr@{}}
		
		\bottomrule
		$\alpha$=0.10\\
		\bottomrule

		$LR$ & 0.993 & 0.920 & 0.593 & 0.193 &  0.087& 0.187 &0.527 &0.880&0.987&0.993\\
		$Z$ & 0.960& 0.900& 0.533&0.167 &0.007 &  0.067& 0.600& 0.96&0.993&1.000 \\
		
		\bottomrule
	\end{tabular}
	
	\caption{Summary of empirical probabilities of concluding the alternative hypothesis for testing $H_{0}: a = a_{0} $ against $H_{1}: a \neq a_{0}$ based on 150 simulation replicates from 1000 independent copies of a true model with parameters $\beta_{MM}=\{\phi_{X}= 0.10, r_{X}=0.25, \phi_{Y}= 0.75, r_{Y}=1.2, a=0.25 \}$, at three significance levels $\alpha \in \{0.01,0.05,0.10\}$. }
	\ra{1.3}
	
		\label{table1}
\end{table}


\begin{table}[H]\centering

	\begin{tabular}{@{}cccccccccccc@{}}
		
		\bottomrule
		$a_{0}$ & 0.10 & 0.25 &0.35 &0.40& 0.45&0.50& 0.55&0.60&0.65 & 0.75  \\
		
		\bottomrule
		$\alpha$=0.01\\
		\bottomrule
		
		$LR$ & 0.993 & 0.967 &0.753  & 0.360 & 0.053 & 0.013 &0.100 &0.560&0.907&0.980\\
		$Z$ &0.967 & 0.960&0.847 &0.393&0.140& 0.000& 0.133&0.653&0.913& 0.960\\
		
	\end{tabular}

	\begin{tabular}{@{}cccccrrrcrrr@{}}
		
		\bottomrule
		$\alpha$=0.05\\
		\bottomrule
		$LR$ & 1.000 & 0.980 &0.933  & 0.627 & 0.227 &  0.033&0.293&0.687&0.960&1.000 \\
		$Z$ & 0.973&0.967 &0.927&0.607 & 0.260& 0.007&0.400&0.807&0.947&0.967\\

	\end{tabular}

	\begin{tabular}{@{}cccccrrrcrrr@{}}
		
		\bottomrule
		$\alpha$=0.10\\
		\bottomrule

		$LR$ &1.000  &0.987  & 0.960 & 0.687&0.260 &0.060 &0.413& 0.813&0.967& 1.000\\
		$Z$ & 0.980&0.967 & 0.960&0.760&0.320 & 0.020&0.540 &0.867&0.960&0.967\\
		
		\bottomrule
	\end{tabular}
	
	\caption{Summary of empirical probabilities of concluding the alternative hypothesis for testing $H_{0}: a = a_{0} $ against $H_{1}: a \neq a_{0}$ based on 150 simulation replicates from 1000 independent copies of a true model with parameters $\beta_{MM}=\{\phi_{X}= 0.10, r_{X}=0.25, \phi_{Y}= 0.75, r_{Y}=1.2, a=0.5 \}$, at three significance levels $\alpha \in \{0.01,0.05,0.10\}$. }
	\ra{1.3}
	
		\label{table2}
\end{table}


\begin{table}[!]\centering
	\ra{1.3}

	\begin{tabular}{@{}cccccccccccc@{}}
		
		\bottomrule
		$a_{0}$& 0.40 &0.50 & 0.60 & 0.65& 0.70 & 0.75 &0.80 & 0.85&0.90 &0.95 \\
		
		\bottomrule
		$\alpha$=0.01\\
		\bottomrule
		
		$LR$ &0.987 & 0.980 &0.767 &0.280  &0.060  & 0.047&0.213&0.593&0.920&0.987 \\
		$Z$ &0.987 & 0.947& 0.740& 0.220&0.047&0.040&0.147&0.367&0.753&0.913\\

	\end{tabular}
	
	\begin{tabular}{@{}cccccrrrcrrr@{}}
		
		\bottomrule
		$\alpha$=0.05\\
		\bottomrule

		$LR$ &0.993  &0.987  &0.900  &0.553  &0.153  &0.087  &0.327&0.780&0.960&1.000 \\
		$Z$ & 0.993& 0.953& 0.907& 0.553& 0.087&0.060 & 0.253&0.620&0.893&0.973\\

	\end{tabular}
	
	\begin{tabular}{@{}cccccrrrcrrr@{}}
		
		\bottomrule
		$\alpha$=0.10\\
		\bottomrule

		$LR$ & 0.993& 0.993 & 0.960 & 0.707 & 0.253 & 0.113 &0.353&0.893&0.980&1.000 \\
		$Z$ & 0.993& 0.960&0.967&0.673&0.187&0.107 &0.307&0.727&0.920 &0.993\\
		
		\bottomrule
	\end{tabular}
	
	\caption{Summary of empirical probabilities of concluding the alternative hypothesis for testing $H_{0}: a = a_{0} $ against $H_{1}: a \neq a_{0}$ based on 150 simulation replicates from 1000 independent copies of a true model with parameters $\beta_{MM}=\{\phi_{X}= 0.10,, r_{X}=0.25, \phi_{Y}= 0.75, r_{Y}=1.2, a=0.75 \}$, at three significance levels $\alpha \in \{0.01,0.05,0.10\}$. }
		\label{table3}
\end{table}

\newpage

\begin{figure} [H]
	\includegraphics[width=0.9\linewidth, height=10cm]{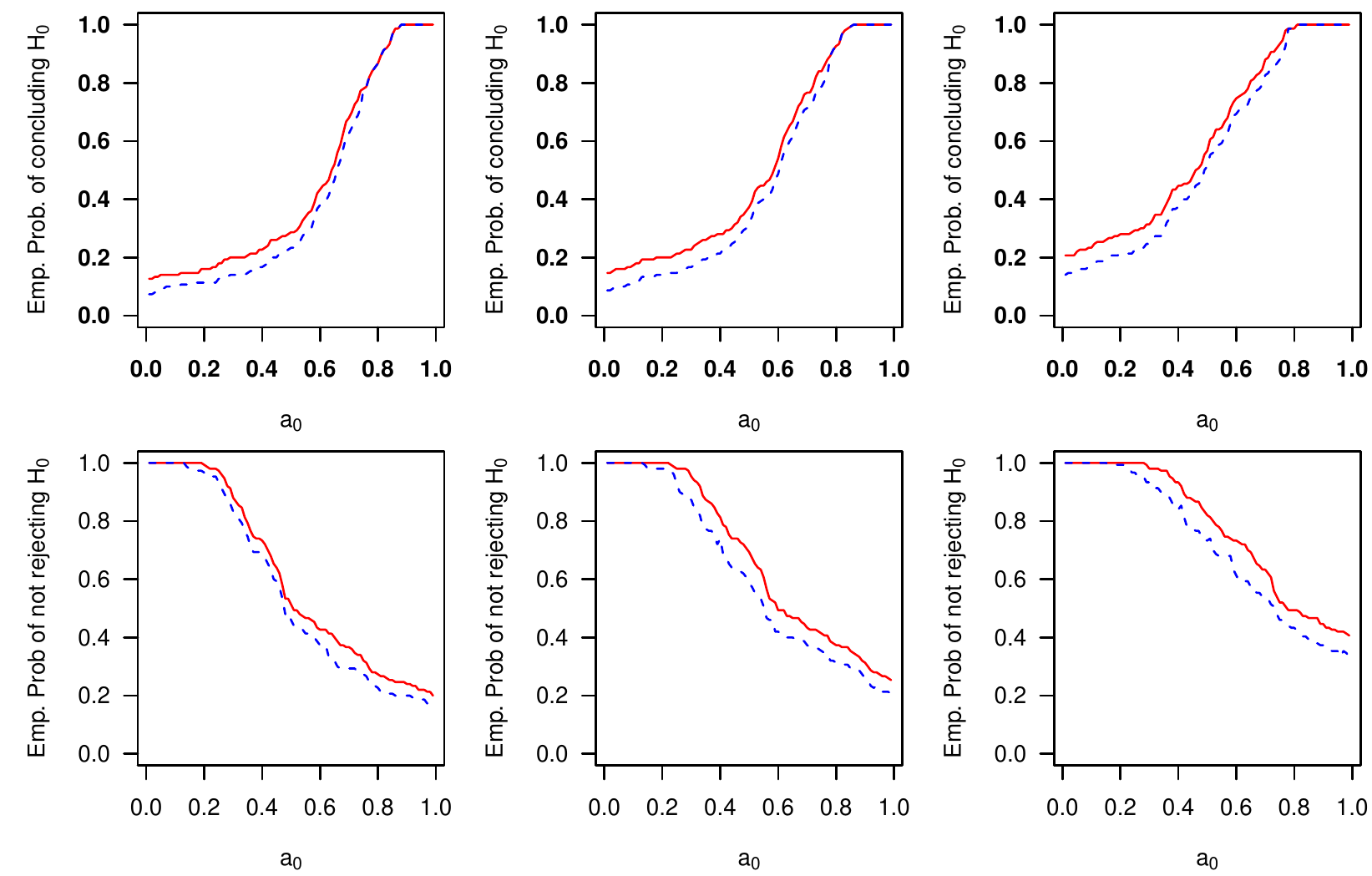}
	
	\caption{ Empirical probabilities of concluding $H_{0}: a = a_{0} $ based on fitting a MM model to 150 simulation replicates from 1000 independent copies from TEG and inverted TEG data. ( Top 
row: MS data with parameters $\beta_{MS}=\{\phi_{X}= 0.10, r_{X}=0.25 \}$,  and bottom row: IMS model with parameters $\beta_{IMS}=\{ \phi_{Y}= 0.75, r_{Y}=1.2 \}$). $50$ locations were randomly 
generated and uniformaly in the square $[0, 2]^{2}$ \vero{(AD case) and $[0\/,1]^2$ in the AI case}, and \vero{$1.000$ realisations} were simulated at the sampled locations. (Red solid line: $Z$ test 
and blue dashed line: $LR$ test). The significance levels from left to right are $\alpha=0.10$, $\alpha=0.05$, and $\alpha=0.01$.}
\label{asymptotic}	
\end{figure}

\vero{We have also explored the AD and AI cases (i.e $a=1$ and $a=0$ respectively).}
The two tests were performed for all $a_0$ values \vero{in the set}  $\{0.01,0.02,...,0.99\}$. \vero{We have computed the probability of concluding $H_0: \ a=a_0$ while the true parameter is $1$ (AD 
case) or $0$ (AI case).} For this purpose a MM model were fitted for each 150 simulation replicates from a MS TEG and inverted TEG processes with parameters $\beta_{MS}=\{ \phi_{X}= 0.10, 
r_{X}=0.25\}$, $\beta_{IMS}=\{r_{Y}=1.2, \phi_{Y}= 0.75\}$, repectively, and a moderately sized data with $K$ = 50 sites and $M$=1000 independent observations. The spatial domains respectively are 
$[0,1]^{2}$ \vero{(AD case)} and $[0,2]^{2}$ \vero{(AI case). The results are presented in Figure \ref{asymptotic} and Tables \ref{table4} and \ref{table5}.}


\begin{table}[!]\centering
	\ra{1.3}
	
	\begin{tabular}{@{}cccccccccccc@{}}
		
		\bottomrule
		$a_{0}$& 0.01 &0.10 & 0.20 & 0.30& 0.40 & 0.50 &0.60 & 0.70&0.80 &0.90 \\
		
		\bottomrule
		$\alpha$=0.01\\
		\bottomrule
			$LR$ &0.140 &0.173 & 0.213&0.247 &0.373& 0.513&0.693&0.820&0.987&1.000\\
		$Z$ &0.207 & 0.233 &0.280&0.313 &0.447  & 0.567&0.747&0.880&0.987&1.000 \\

	\end{tabular}
	
	\begin{tabular}{@{}cccccrrrcrrr@{}}
		
		\bottomrule
		$\alpha$=0.05\\
		\bottomrule

		$LR$ & 0.087 & 0.113 & 0.140 & 0.167 & 0.213 & 0.320&0.487&0.713&0.927&1.000 \\
		$Z$ & 0.147& 0.173& 0.200& 0.227& 0.280&0.373 & 0.540&0.767&0.927&1.000\\

	\end{tabular}
	
	\begin{tabular}{@{}cccccrrrcrrr@{}}
		
		\bottomrule
		$\alpha$=0.10\\
		\bottomrule

		$LR$ &0.073 & 0.100 & 0.113 &0.140 &  0.167& 0.233 &0.380&0.613&0.867&1.000 \\
		$Z$ &0.127 & 0.140&0.160&0.200&0.227&0.287 &0.433&0.680&0.867 &1.000\\
		
		\bottomrule
	\end{tabular}
	
	\caption{Summary of empirical probabilities of concluding the null hypothesis for testing $H_{0}: a = a_{0} $ against $H_{1}: a \neq a_{0}$ based on 150 simulation replicates from 1000 independent copies of a true MS model with parameters $\beta_{MS}=\{\phi_{X}= 0.10, r_{X}=0.25 \}$, at three significance levels $\alpha \in \{0.01,0.05,0.10\}$. }
	\label{table4}
\end{table}


\begin{table}[!]\centering
	\ra{1.3}
	
	\begin{tabular}{@{}cccccccccccc@{}}
		
		\bottomrule
		$a_{0}$&0.10 & 0.20 & 0.30& 0.40 & 0.50 &0.60 & 0.70&0.80 &0.90&0.99 \\
		
		\bottomrule
		$\alpha$=0.01\\
		\bottomrule
		$LR$ & 1.000&0.993 & 0.933&0.840 &0.733& 0.613&0.540&0.433&0.373&0.333\\
		$Z$ &1.000 & 1.000 &0.980&0.933 &0.820 & 0.733&0.633&0.493&0.447&0.407 \\

	\end{tabular}
	
	\begin{tabular}{@{}cccccrrrcrrr@{}}
		
		\bottomrule
		$\alpha$=0.05\\
		\bottomrule

		$LR$ &1.000  & 0.980 &0.873& 0.740 & 0.587 &0.420 & 0.373& 0.307& 0.267& 0.200\\
		$Z$ & 1.000& 1.000& 0.953& 0.813& 0.693&0.493 & 0.433&0.373&0.313&0.253\\

	\end{tabular}
	
	\begin{tabular}{@{}cccccrrrcrrr@{}}
		
		\bottomrule
		$\alpha$=0.10\\
		\bottomrule

		$LR$ &1.000&0.967  &  0.833& 0.680&0.487  & 0.367 &0.293& 0.220&0.200&0.167\\
		$Z$ &1.000 & 0.993&0.880&0.733&0.513&0.427 &0.367&0.273&0.240&0.200\\
		
		\bottomrule
	\end{tabular}
	
	\caption{Summary of empirical probabilities of concluding the null hypothesis for testing $H_{0}: a = a_{0} $ against $H_{1}: a \neq a_{0}$ based on 150 simulation replicates from 1000 independent copies of a true IMS model with parameters $\beta_{IMS}=\{ \phi_{Y}= 0.75, r_{Y}=1.2 \}$, at three significance levels $\alpha \in \{0.01,0.05,0.10\}$. }
	\label{table5}
\end{table}

Despite the interesting results, that show the increase of the empirical probability of concluding $H_{0}$, \vero{both tests fail to identify precisely AD or AI. They could be used as indication 
of AD or AI but are not sensitive.}


\section{Rainfall Data Example}\label{sec:rainfall}

\begin{figure} [H]
	
	\includegraphics[width=0.9\linewidth, height=10cm]{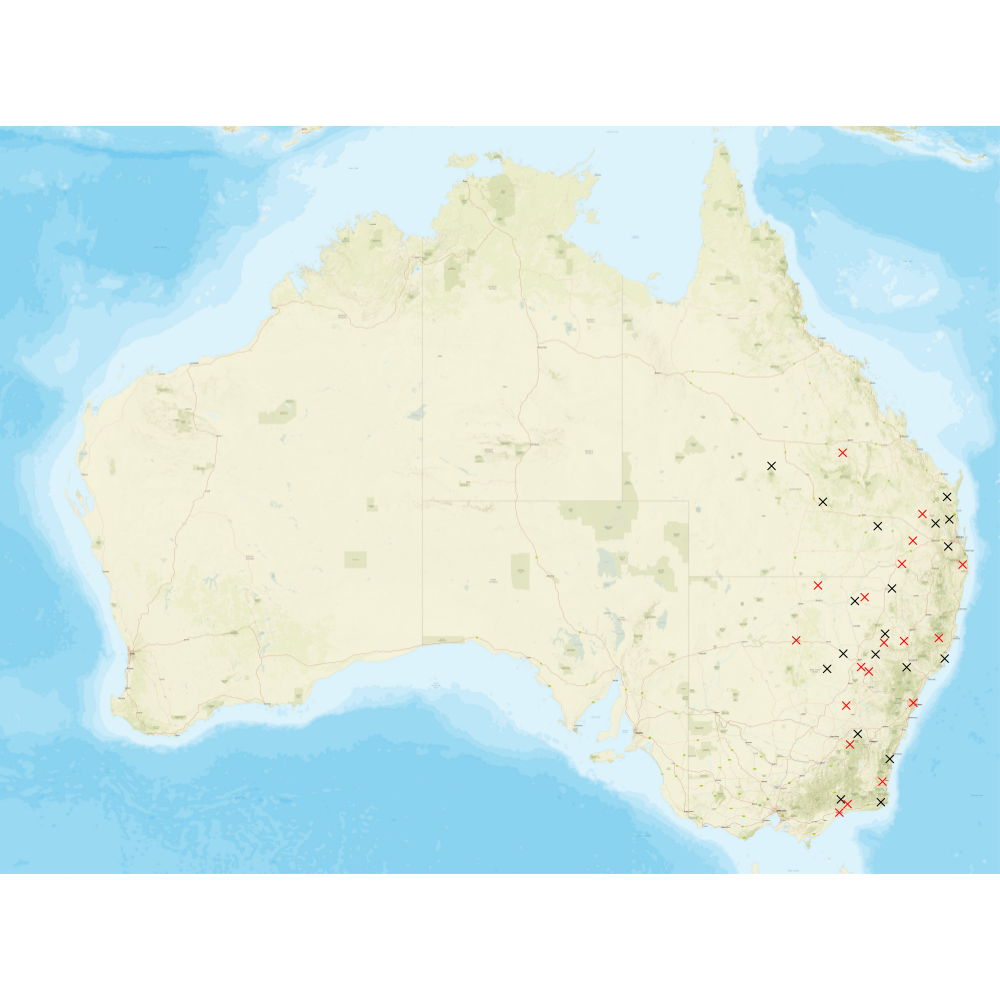}
	
	\caption{Geographical locations of the 38 meteorological stations located in the East of Australia.  \textcolor{black} {Red crosses represent stations in  group A, while black ones denote group B}.}
	
	\label{Map}
	
\end{figure}

The data analysed in this section are daily rainfall  amounts in (millimetres) over the years 1972-2016 occurring during April-September at 38 sites in the East of Australia whose locations are shown 
in Fig. \ref{Map}. The altitude of the sites varying from 4 to 552 meters above mean sea level. The sites are separated by distances \vero{from $34$km to $1383$km}. This data set is freely avaliable 
on http://www.bom.gov.au/climate/data/. \\

\begin{figure} [H]
	\includegraphics[width=0.99\linewidth]{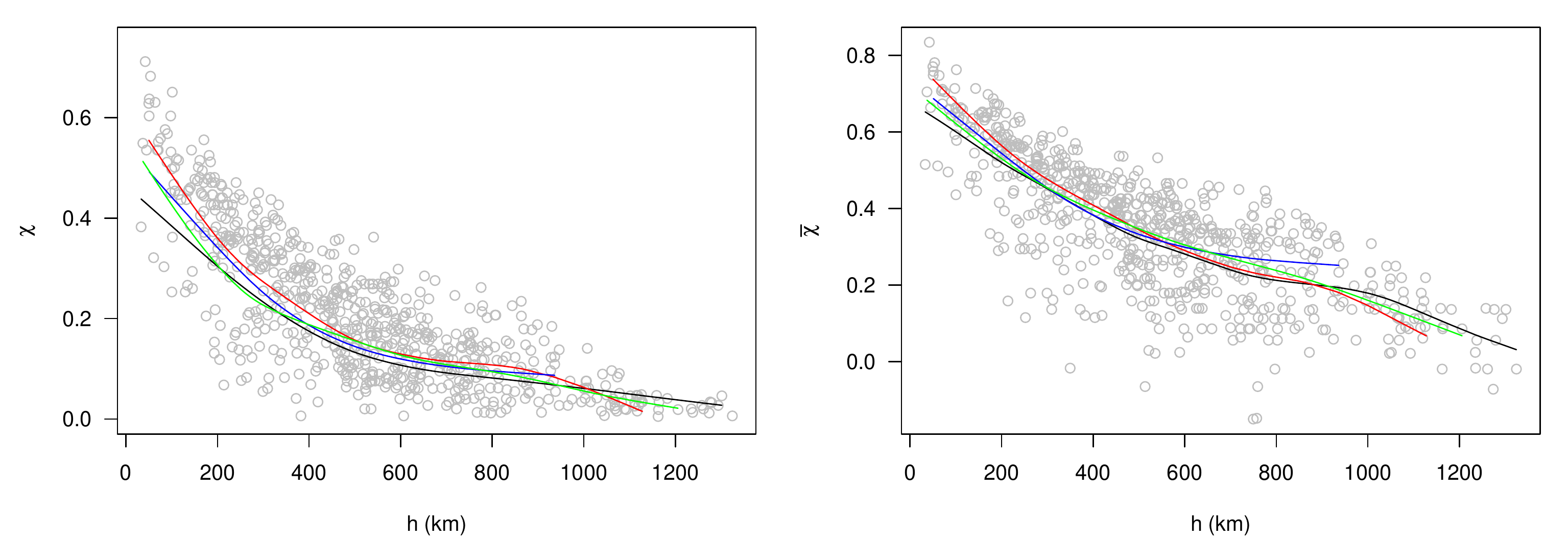}
	
	\caption{ Pairwise empirical estimates of $\chi$ (left panel) and $\bar{\chi}$ (right panel) versus distance at threshold $u = 0.970$. Grey points are empirical pairwise 
estimates for all data pairs. Colored lines are the loess smoothed values of the empirical estimates in different directional sectors: black line $(- \pi/8,  \pi/8]$, red line $( \pi/8,  3\pi/8]$, 
blue line $( 3\pi/8,  5\pi/8]$, and green line $( 5\pi/8,  7\pi/8]$. }
	\label{Isotropy}	
\end{figure}

Following the approach of \cite{bacro2016flexible}, \vero{we have done a graphical exploration using the coefficients $\chi$ and $\bar{\chi}$, in order to evaluate a} possible anisotropy \vero{in the 
data.  Figure \ref{Isotropy} is} based on the empirical estimates of the functions $\chi(h,u)$ and $\bar{\chi}{(h,u)}$ in different directional sectors $(- \pi/8,  \pi/8]$,  $( \pi/8,  3\pi/8]$, 
$( 3\pi/8,  5\pi/8]$, and $( 5\pi/8,  7\pi/8]$, where $0$ represents the northing direction. \vero{We can conclude that there is no evidence for anisotropy and we shall consider that the data comes 
from an isotropic process}.

A way to apply our testing approach is based on dividing the daily rainfall dataset from the 38 sites into two groups A and B (see \vero{Figure} \ref{Map}). \vero{We shall consider} five 
models \vero{which belong} to the three classes: MM, MS, and IMS. \vero{These models are first fitted on data from} group A. The composite likelihood information criterion (CLIC) 
\cite{varin2005note}, \vero{defined} as CLIC= $-2 [p\ell(\hat\vartheta) - tr\{\mathcal{J}(\hat\vartheta)\mathcal{H}^{-1}(\hat\vartheta)\}]$ is used to \vero{choose} the best \vero{fitted} model. 
Lower values of CLIC indicate \vero{a} better fit. Then \vero{we apply our test with the best fitted} model for group A. \vero{Let $a_0$ be the mixing parameter of the best fitted model for group A 
(we found $a_0=0.23$). Our test will be done on the group B data, and will be $H_0: a=a_0$ vs $H_1: a\neq a_0$.} \\

 The fitted models are:\\

$ \mathbf{M_{1}}$: a MM model where $X$ is a TEG process with an exponential correlation function 
$ \rho(h) = \exp (- \lVert {h} \rVert / \phi_{X})$, $\phi_{X}>0$.  $\mathcal{A}_{X}$ \vero{is a disk of fixed} and unknown radius $r_{X}$, and $Y$ is an inverted TEG process with exponential 
correlation function $ \rho(h) = \exp (- \lVert {h} \rVert / \phi_{Y})$, $\phi_{Y}>0$, and $\mathcal{A}_{Y}$ is a disk with fixed and unknown radius $r_{Y}$.\\

$\mathbf{M_{2}}$: a MM model where $X$ is a TEG process as in $ \mathbf{M_{1}}$. $Y$ is an isotropic inverted Smith process where $\Sigma$ is a diagonal matrix ($\sigma_{12}=0$) with $\sigma_{11}^{2}=\sigma_{22}^{2}=\phi_{Y}^{2}$, i.e., $\gamma(h)= ( \lVert {h} \rVert / \phi_{Y})$.\\
 
$\mathbf{M_{3}}$: \vero{a MS} TEG process described as $X$ in $ \mathbf{M_{1}}$.\\

$\mathbf{M_{4}}$: a MS isotropic Smith process where $\Sigma$ is a diagonal matrix ($\sigma_{12}=0$) with $\sigma_{11}^{2}=\sigma_{22}^{2}=\phi_{X}^{2}$, i.e., $\gamma(h)= ( \lVert {h} \rVert / \phi_{X})$.\\ 
 
$\mathbf{M_{5}}$: the inverted Smith process described as $Y$ in $ \mathbf{M_{2}}$.\\
 
\vero{The considered models have unit Fréchet marginal distributions. We fit a GEV distribution on each site and then transform the marginal laws to unit Fréchet using}
$$\textcolor{black}{\hat{F}(z)= \left(1+\hat{\xi} \frac{z-\hat{\mu}}{\hat{\sigma}}\right)^{1/\hat{\xi}}\/},$$ 
where $(\hat{\xi},\hat{\mu},\hat{\sigma})$ are the estimated parameters of the GEV distribution. The censored pairwise likelihood approach (\ref{pairwise1}) is used \vero{in order to estimate the 
parameters. The threshold is $u=-1/\log(0.90)$ and equal weights are used. The matrices} $ \mathcal{H}(\vartheta)$ and $\mathcal{J}(\vartheta)$ and the related quantities CLIC and standard errors 
are obtained by Monte Carlo procedure through simulating data with $M =1500$ independent draws at the sampled 19 sites from the fitted model.\\

Our results are \vero{summarised in Table \ref{table 6}}. The best-fitting model for group A, as judged by CLIC, is the hybrid dependence model $M_{2}$.  \textcolor{black} {The mixing coefficient for the other hybrid 
model $M_{1}$ is very close to one, which indicates that there is no mixture between the max-stable process and the asymptotically independent one, so the asymptotic independence components are not 
identifiable. This fact affects the values of the estimates. Moreover, model $M_{1}$ reduces to model $M_{3}$}.
 
 
\begin{table}[H]\centering
	\ra{1.3}
	
	\begin{tabular}{@{}l l l l l l l@{}}
		
		\bottomrule
		Model &${\hat{\phi}}_{X}$  & ${\hat{r}}_{X}$& $\hat{{a}}$ & ${\hat{\phi}}_{Y}$ &${\hat{r}}_{Y}$& CLIC  \\
		
		\bottomrule
	    $M_{1}$ & 302.09 (67.79)& 753.67 (191.02) &  $\simeq1$(0.01) &2064.71 (542.84) &970.04 (188.28)& 1951654\\
	  
		$M_{2}$ & 43.77 (31.84)& 94.81 (45.96) & 0.23 (0.06) &1111.53 (430.98)&- & ${1950330}^{*}$\\ 
		 $M_{3}$ &303.09 (68.52)& 751.44 (189.10) &-&-&-&1951654 \\              
		$M_{4}$ & 130.05 (36.82) & -&-&- &-&1968505\\
		$M_{5}$ &-& - &-&337.40 (86.41)  &-& 1963187\\
			\bottomrule

	\end{tabular}
	\caption{Parameter estimates of selected dependence models fitted to the daily rainfall data from group A, along with the composite likelihood criterion (CLIC) and standard errors reported between parentheses. }
	
\label{table 6}
\end{table}

\vero{For data from group B, we have considered the two statistics $Z$ and $LR$ described in this paper,} \textcolor{black} {to test if the hybrid model $M_{2}$ can be used to make inference for this group, i.e., testing $H_{0}: a=0.23$ versus $H_{1}: a\neq 0.23$. We obtain the calculated values for $Z$ and $LR$, $Z=0.73$ ($p-$value$=0.47$), ${LR=7.72}$, $\lambda=14.07$, ($p-$value$={0.458}$)}.
\\

\vero{This leads us to retain $H_0$ and thus that there is no differences in the mixing parameter between the two groups. We have also performed an independent two-samples $Z$-test: let $a_A$ (resp. 
$a_B$ be the mixing parameter for group A (resp. group B), consider $H_{0}: a_{A} = a_{B}$. The
statistic 
$$Z_C=\frac{\hat{a}_{A}-\hat{a}_{B}}{ SE_{\hat{a}_{A}-\hat{a}_{B}}} $$
where $SE$ stands for the estimated standard error. The calculated value of $Z_C$ is $0.17$ with $p$-value $= 0.86$, leading to retain $H_0$ and conclude that there is no significant difference 
between  the two mixing coefficient.} 
\\
\vero{Nevertheless, these conclusions are subordinated to the assumption that both groups A and B have the same underlying model M2.}

\section{Conclusion}\label{sec:conclusion}

In this paper we have considered hypothesis testing for the mixing coefficient of a MM models \vero{proposed in} \cite{doi:10.1093/biomet/asr080} using two statistics the $Z$ and the $LR$ when a 
censored pairwise likelihood is \vero{employed} for inferential purposes.\\

The two statistics has emerged as an efficient tools for testing hypothesis \vero{on} the mixing coefficient with better performance \vero{achieved} by $LR$, but \vero{with} the drawback of a 
nonstandard asymptotic distribution at the boundaries, since the number of nuisance parameters is different between the two hypothesis \vero{and it also requires heavier computations. Our procedure seems to be a performant validation tool.}  \\

\vero{One other} drawback \vero{of our work is that the proposed tests model-dependent}. In the future, we plan to propose a free-model test, \vero{using the F-}madogram. 
\\
\ \\ 
\ \\
\underline{Acknowledgments.} We are grateful to Jean-Noël Bacro, Carlo Gaetan and Gwladys Toulemonde for giving their estimations C codes which we used as a base for computing the statistics of our tests.  This work was supported by the LABEX MILYON
(ANR-10-LABX-0070) of Université de Lyon, within the program ”Investissements
d’Avenir” (ANR-11-IDEX-0007) operated by the French National Research
Agency (ANR). It was also supported by the CERISE LEFE-INSU projcet. 
\bibliographystyle{plain}
\bibliography{Mybib}

\begin{thebibliography}{10}

\bibitem{bacro2010testing}
Jean-No{\"e}l Bacro, Liliane Bel, and Christian Lantu{\'e}joul.
\newblock Testing the independence of maxima: from bivariate vectors to spatial
  extreme fields.
\newblock {\em Extremes}, 13(2):155--175, 2010.

\bibitem{bacro2016flexible}
Jean-Noel Bacro, Carlo Gaetan, and Gwladys Toulemonde.
\newblock A flexible dependence model for spatial extremes.
\newblock {\em Journal of Statistical Planning and Inference}, 172:36--52,
  2016.

\bibitem{cattelan2016empirical}
Manuela Cattelan and Nicola Sartori.
\newblock Empirical and simulated adjustments of composite likelihood ratio
  statistics.
\newblock {\em Journal of Statistical Computation and Simulation},
  86(5):1056--1067, 2016.

\bibitem{chandler2007inference}
Richard~E Chandler and Steven Bate.
\newblock Inference for clustered data using the independence loglikelihood.
\newblock {\em Biometrika}, pages 167--183, 2007.

\bibitem{coles1999dependence}
Stuart Coles, Janet Heffernan, and Jonathan Tawn.
\newblock Dependence measures for extreme value analyses.
\newblock {\em Extremes}, 2(4):339--365, 1999.

\bibitem{cooley2006variograms}
Dan Cooley, Philippe Naveau, and Paul Poncet.
\newblock Variograms for spatial max-stable random fields.
\newblock {\em Dependence in probability and statistics}, pages 373--390, 2006.

\bibitem{davies1977hypothesis}
Robert~B Davies.
\newblock Hypothesis testing when a nuisance parameter is present only under
  the alternative.
\newblock {\em Biometrika}, 64(2):247--254, 1977.

\bibitem{davies1987hypothesis}
Robert~B Davies.
\newblock Hypothesis testing when a nuisance parameter is present only under
  the alternatives.
\newblock {\em Biometrika}, pages 33--43, 1987.

\bibitem{davison2012geostatistics}
AC~Davison and MM~Gholamrezaee.
\newblock Geostatistics of extremes.
\newblock In {\em Proc. R. Soc. A}, volume 468, pages 581--608, 2012.

\bibitem{davison2013geostatistics}
Anthony~C Davison, Rapha{\"e}l Huser, and Emeric Thibaud.
\newblock Geostatistics of dependent and asymptotically independent extremes.
\newblock {\em Mathematical Geosciences}, 45(5):511--529, 2013.

\bibitem{davison2012statistical}
Anthony~C Davison, Simone~A Padoan, and Mathieu Ribatet.
\newblock Statistical modeling of spatial extremes.
\newblock {\em Statistical science}, pages 161--186, 2012.

\bibitem{de2006spatial}
L.~De~Haan and T.~T. Pereira.
\newblock Spatial extremes: Models for the stationary case.
\newblock {\em The annals of statistics}, pages 146--168, 2006.

\bibitem{de1984spectral}
Laurens De~Haan.
\newblock A spectral representation for max-stable processes.
\newblock {\em The annals of probability}, pages 1194--1204, 1984.

\bibitem{geys1999pseudolikelihood}
Helena Geys, Geert Molenberghs, and Louise~M Ryan.
\newblock Pseudolikelihood modeling of multivariate outcomes in developmental
  toxicology.
\newblock {\em Journal of the American Statistical Association},
  94(447):734--745, 1999.

\bibitem{huser2014space}
Rapha{\"e}l Huser and AC~Davison.
\newblock Space--time modelling of extreme events.
\newblock {\em Journal of the Royal Statistical Society: Series B (Statistical
  Methodology)}, 76(2):439--461, 2014.

\bibitem{joe1993parametric}
Harry Joe.
\newblock Parametric families of multivariate distributions with given margins.
\newblock {\em Journal of multivariate analysis}, 46(2):262--282, 1993.

\bibitem{kabluchko2009stationary}
Zakhar Kabluchko, Martin Schlather, and Laurens De~Haan.
\newblock Stationary max-stable fields associated to negative definite
  functions.
\newblock {\em The Annals of Probability}, pages 2042--2065, 2009.

\bibitem{kent1982robust}
John~T Kent.
\newblock Robust properties of likelihood ratio tests.
\newblock {\em Biometrika}, 69(1):19--27, 1982.

\bibitem{kuonen1999miscellanea}
Diego Kuonen.
\newblock Miscellanea. saddlepoint approximations for distributions of
  quadratic forms in normal variables.
\newblock {\em Biometrika}, 86(4):929--935, 1999.

\bibitem{ledford1996statistics}
Anthony~W Ledford and Jonathan~A Tawn.
\newblock Statistics for near independence in multivariate extreme values.
\newblock {\em Biometrika}, 83(1):169--187, 1996.

\bibitem{lindsay1988composite}
Bruce~G Lindsay.
\newblock Composite likelihood methods.
\newblock {\em Contemporary mathematics}, 80(1):221--39, 1988.

\bibitem{opitz2013extremal}
Thomas Opitz.
\newblock Extremal t processes: Elliptical domain of attraction and a spectral
  representation.
\newblock {\em Journal of Multivariate Analysis}, 122:409--413, 2013.

\bibitem{pace2011adjusting}
Luigi Pace, Alessandra Salvan, and Nicola Sartori.
\newblock Adjusting composite likelihood ratio statistics.
\newblock {\em Statistica Sinica}, pages 129--148, 2011.

\bibitem{padoan2010likelihood}
Simone~A Padoan, Mathieu Ribatet, and Scott~A Sisson.
\newblock Likelihood-based inference for max-stable processes.
\newblock {\em Journal of the American Statistical Association},
  105(489):263--277, 2010.

\bibitem{rotnitzky1990hypothesis}
Andrea Rotnitzky and Nicholas~P Jewell.
\newblock Hypothesis testing of regression parameters in semiparametric
  generalized linear models for cluster correlated data.
\newblock {\em Biometrika}, pages 485--497, 1990.

\bibitem{schlather2002models}
Martin Schlather.
\newblock Models for stationary max-stable random fields.
\newblock {\em Extremes}, 5(1):33--44, 2002.

\bibitem{smith1990max}
Richard~L Smith.
\newblock Max-stable processes and spatial extremes.
\newblock {\em Unpublished manuscript}, 205, 1990.

\bibitem{thibaud2013threshold}
E~Thibaud, R~Mutzner, and AC~Davison.
\newblock Threshold modeling of extreme spatial rainfall.
\newblock {\em Water resources research}, 49(8):4633--4644, 2013.

\bibitem{varin2008composite}
Cristiano Varin.
\newblock On composite marginal likelihoods.
\newblock {\em AStA Advances in Statistical Analysis}, 92(1):1--28, 2008.

\bibitem{varin2011overview}
Cristiano Varin, Nancy Reid, and David Firth.
\newblock An overview of composite likelihood methods.
\newblock {\em Statistica Sinica}, pages 5--42, 2011.

\bibitem{varin2005note}
Cristiano Varin and Paolo Vidoni.
\newblock A note on composite likelihood inference and model selection.
\newblock {\em Biometrika}, pages 519--528, 2005.

\bibitem{doi:10.1093/biomet/asr080}
Jennifer~L. Wadsworth and Jonathan~A. Tawn.
\newblock Dependence modelling for spatial extremes.
\newblock {\em Biometrika}, 99(2):253, 2012.

\end{thebibliography}

\end{document}